\documentclass[a4paper,10pt]{amsart}
\usepackage{latexsym}
\usepackage{amssymb}
\usepackage{amsfonts}
\usepackage{amsmath}
\usepackage{indentfirst}
\usepackage{graphicx}
\usepackage{epsfig}

\usepackage{epic, eepic}

\theoremstyle{plain}
\newtheorem{teor}{Theorem}[section]
\newtheorem{cor}[teor]{Corollary}

\newtheorem{lemma}[teor]{Lemma}
\newtheorem{prop}[teor]{Proposition}
\theoremstyle{definition}
\newtheorem{defn}[teor]{Definition}

\theoremstyle{remark}

\newcommand\sig{\ensuremath{\sigma}}
\newcommand\taup{\ensuremath{{\tau'}}}
\newcommand\ptau{\ensuremath{{'\kern-.4ex\tau}}}
\newcommand\ptp{\ensuremath{{'\kern-.4ex\tau'}}}

\newcommand\x{\ensuremath{{\mathcal C}}}
\newcommand\pC{\ensuremath{{'\kern-.2ex C}}}
\newcommand\Cp{\ensuremath{C}'}

\newcommand\mf{M\"obius function}

\def\cls{{\mathcal S}}
\def\mp{{\mathcal P}}

\newcommand\ist{\ensuremath{[\sigma,\tau]}}

\newcommand\mst{\ensuremath{\mu(\sigma,\tau)}}
\newcommand\msx{\ensuremath{\mu(\sigma,\x)}}
\def\m#1,#2,{\ensuremath{\mu(#1,#2)}}
\def\mb#1,#2,{\ensuremath{\bar\mu(#1,#2)}}

\author[A. Bernini, L. Ferrari and E. Steingr\'imsson]{Antonio
  Bernini, Luca Ferrari, and Einar Steingr\'imsson}

\address{Dipartimento di Sistemi e Informatica, University of Firenze, Italy}
\email{bernini@dsi.unifi.it}

\address{Dipartimento di Sistemi e Informatica, University of Firenze, Italy}
\email{ferrari@dsi.unifi.it}

\address{Department of Computer and Information Sciences, University
  of Strathclyde, Glasgow G1 1XH, UK}
\email{einar.steingrimsson@cis.strath.ac.uk}

\thanks{Steingr\'imsson was supported by grant no.\ 090038012 from the
  Icelandic Research Fund.}

\title[The M\"{o}bius function of the consecutive pattern poset]{The
  M\"{o}bius function of the\\ consecutive pattern poset}

\usepackage[usenames]{color}
\usepackage{colordvi}

\begin{document}

\begin{abstract}
  An occurrence of a consecutive permutation pattern $p$ in a
  permutation $\pi$ is a segment of consecutive letters of $\pi$ whose
  values appear in the same order of size as the letters in $p$.  The
  set of all permutations forms a poset with respect to such pattern
  containment.  We compute the M\"obius function of intervals in this
  poset, providing what may be called a complete solution to the
  problem.  For most intervals our results give an immediate answer to
  the question.  In the remaining cases, we give a polynomial time
  algorithm to compute the M\"obius function.  In particular, we show
  that the M\"obius function only takes the values $-1$, $0$ and ~$1$.
\end{abstract}

\maketitle
\thispagestyle{empty}

\section{Preliminaries and introduction}

For the poset of classical permutation patterns, the first results
about its M\"obius function were obtained in \cite{SV}.  Further
results appear in \cite{ST} and \cite{BJJS}.  The general problem in
this case of classical patterns seems quite hard.  In contrast, the
poset of consecutive pattern containment has a much simpler structure.
In this paper we compute the M\"obius function of that poset and
provide what may be called a complete solution to the problem.  In
most cases our results give an immediate answer.  In the remaining
cases, we give a polynomial time recursive algorithm to compute the
M\"obius function.  In particular, we show that the M\"obius function
only takes the values $-1$, $0$ and ~$1$.

Unless otherwise specified, all permutations in this paper are taken
to be of the set $[d]=\{1,2,\ldots,d\}$ for some positive integer $d$.
We denote by $\cls_d$ the set of all such permutations for a given
$d$.  An occurrence of a \emph{consecutive pattern} $\sig=a_1a_2\dots
a_k$ in a permutation $\tau=b_1b_2\dots b_n$ is a subsequence
$b_{i}b_{i+1}\dots b_{i+k-1}$ in $\tau$, whose letters appear in the
same order of size as the letters in $\sig$.  As an example, there are
three occurrences of the consecutive pattern $231$ in the permutation
563724891, namely 563, 372 and 891.  On the other hand, the
permutation 253641 \emph{avoids} 231, since it contains no occurrence
of that pattern.

Consecutive permutation patterns are special cases of the generalized
permutation patterns introduced in \cite{BS}, and they are not to be
confused with the classical permutation patterns, whose occurrences in
a permutation do not have to be contiguous.  The enumerative
properties of occurrences of various consecutive permutation patterns
were first studied systematically in \cite{EN}, but these results will
not concern us, as there seems to be no connection between them and
the M\"obius function studied here.

The set of all permutations forms a poset $\mathcal P$ with respect to
consecutive pattern containment.  In other words, if $\sig\in\cls_k$
and $\tau\in\cls_n$, then $\sig\leq\tau$ in $\mathcal P$ if $\sig$
occurs as a consecutive pattern in $\tau$.  We write $\sig<\tau$ if
$\sig\le\tau$ and $\sig\neq\tau$.  As usual in poset terminology, a
permutation $\tau$ \emph{covers} $\sig$ (and \emph{$\sig$ is covered
  by $\tau$}) if $\sig<\tau$ and there is no permutation $\pi$ such
that $\sig<\pi<\tau$.  Note that if $\tau$ covers $\sig$ then
$|\tau|-|\sig|=1$, where $|\pi|$ is the length of $\pi$. The interval
$[x,y]$ in a poset $P$, where $x$ and $y$ are elements of
$P$, is defined by $[x,y]=\{z\in P\;|\; x\le z\le y\}$. The
\emph{rank} of an interval $\ist$ in $\mp$ is the difference
$|\tau|-|\sig|$. The rank of an element $\pi$ in $\ist$ is defined to
be the rank of the interval $[\sig,\pi]$.

A \emph{filter} in a poset $P$ is a set $S\subseteq P$ such that if
$x>y$ and $y\in S$, then $x\in S$.  An \emph{ideal} is a set
$S\subseteq P$ such that if $x<y$ and $y\in S$, then $x\in S$.  A
\emph{principal filter} is a filter with a single minimal element, and
a \emph{principal ideal} is an ideal with a single maximal  element.
In each case, the single minimal/maximal element is said to
\emph{generate} the filter/ideal.

In the poset $\mathcal P$ consider the interval $[\sigma,\tau]$.  Our aim is to
compute $\mu(\sigma,\tau)$, where $\mu$ is the M\"{o}bius function of
the incidence algebra of $\mathcal P$.  The M\"obius function is
recursively defined by setting $\mu(x,x)=1$ for all $x$, and, if
$x\neq y$,
\begin{equation}\label{def-mob}
\mu(x,y)= -\sum_{x\le z<y}{\mu(x,z)}.
\end{equation}
In particular, if $x\not\leq y$, then $\mu(x,y)=0$.

It is well known (and follows easily from Philip Hall's Theorem,
see \cite[Proposition 3.8.5]{St}) that the M\"obius function of an
interval $I$ is equal to the M\"obius function of the interval
$I'$ obtained by ``turning $I$ upside-down,'' that is, the
interval $I'$ for which $x$ is declared to be less than or equal
to $y$ when $y\le x$ in $I$. Thus, we also have
\begin{equation}\label{up-down-mob}
\mu(x,y)= -\sum_{x<z\le y}{\mu(z,y)}.
\end{equation}

Examples of the two ways to compute the M\"obius function, described
in Equations~(\ref{def-mob}) and~(\ref{up-down-mob}), are given in
Figure~\ref{fig-updown}. Denoting with $[x,y]$ the interval whose
Hasse diagram is depicted in the figure, the label of an element $z$
on the left is $\mu(x,z)$, whereas the label of $z$ on the right is
$\mu(z,y)$.

Note also that it follows from either~(\ref{def-mob})
or~(\ref{up-down-mob}) that the sum of $\mu(x,z)$ over \emph{all}
elements $z$ in an interval $[x,y]$ is 0, that is
\begin{equation}\label{mob-zero}
\sum_{x\le z\le y}{\mu(x,z)} = 0.
\end{equation}
This identity will be used frequently.

\def\p{\circle*{1.5}}
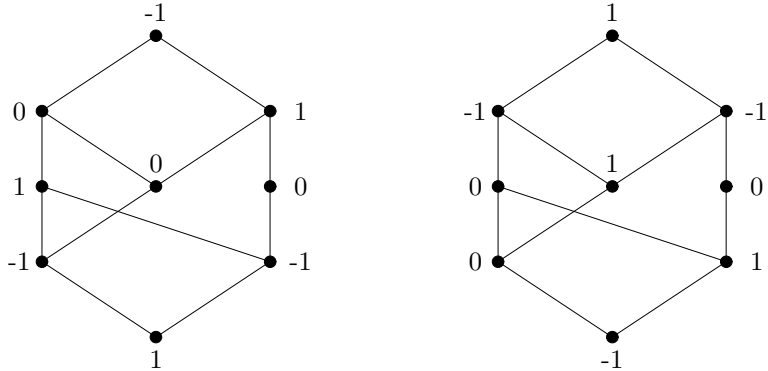
\begin{figure}[htbp]
\begin{center}
\setlength{\unitlength}{1mm}
\begin{picture}(95,50)(-3,-5)
\put(15,40){\p}

\put(0,30){\p}
\put(30,30){\p}

\put(0,20){\p}
\put(15,20){\p}
\put(30,20){\p}

\put(0,10){\p}
\put(30,10){\p}

\put(15,0){\p}

\path(15,0)(0,10)(0,30)(15,40)(30,30)(30,10)(15,0)
\path(0,10)(30,30)
\path(0,20)(30,10)
\path(0,30)(15,20)

\put(15,-3){\makebox(0,0){1}}
\put(-3,10){\makebox(0,0){-1}}
\put(34,10){\makebox(0,0){-1}}

\put(-3,20){\makebox(0,0){1}}
\put(34,20){\makebox(0,0){0}}
\put(15,23){\makebox(0,0){0}}

\put(-3,30){\makebox(0,0){0}}
\put(34,30){\makebox(0,0){1}}

\put(15,43){\makebox(0,0){-1}}

\put(60,0){
\put(15,40){\p}

\put(0,30){\p}
\put(30,30){\p}

\put(0,20){\p}
\put(15,20){\p}
\put(30,20){\p}

\put(0,10){\p}
\put(30,10){\p}

\put(15,0){\p}

\path(15,0)(0,10)(0,30)(15,40)(30,30)(30,10)(15,0)
\path(0,10)(30,30)
\path(0,20)(30,10)
\path(0,30)(15,20)

\put(15,-3){\makebox(0,0){-1}}

\put(-3,10){\makebox(0,0){0}}
\put(34,10){\makebox(0,0){1}}

\put(-3,20){\makebox(0,0){0}}
\put(34,20){\makebox(0,0){0}}
\put(15,23){\makebox(0,0){1}}

\put(-3,30){\makebox(0,0){-1}}
\put(34,30){\makebox(0,0){-1}}

\put(15,43){\makebox(0,0){1}}
}
\end{picture}
\caption{\label{fig-updown} Computing the M\"obius function of an
  interval from bottom to top (left) and from top to bottom (right).}
\end{center}
\end{figure}

\begin{defn}
  Given a sequence of distinct integers $s=s_1s_2\dots s_d$, \emph{the
    standard form of $s$} is the permutation $\pi=a_1a_2\dots a_d$ of
  $\{1,2,\dots,d\}$ that is \emph{order isomorphic} to $s$, that is,
  whose letters appear in the same order of size as those of $s$.
\end{defn}
For example, the standard form of both 4731 and 6842 is 3421, so all
three are order isomorphic.  Note that, since an occurrence of a
consecutive pattern $\rho$ in a permutation $\pi$ only depends on the
relative order of the letters in $\rho$, we may as well talk about the
occurrence of an arbitrary sequence of distinct integers as a
consecutive pattern in $\pi$.  For example, the subsequence 352 in
435261 is an occurrence of the sequence 451.

\begin{defn}
Suppose $\sig$ occurs in $\tau=a_1a_2\dots a_n$.  If $a_{i+1}$ is the
leftmost letter of $\tau$ involved in any occurrence of $\sig$ in
$\tau$, we say that $\tau$ has a \emph{left tail of length $i$ with
  respect to $\sig$}.  Analogously, $\tau$ has a \emph{right tail of
  length $j$ with respect to $\sig$} if $a_{n-j}$ is the rightmost
letter of $\tau$ involved in any occurrence of $\sig$ in $\tau$.  If
it is clear from the context what $\sig$ is, we simply talk about left
and right tails of $\tau$.
\end{defn}

The following definition is borrowed from the theory of codes.

\begin{defn}\label{def_pre_patt}
  Given a permutation $\tau$, its \emph{prefix} (resp.\ \emph{suffix})
  \emph{pattern of length}~$k$ is the permutation of length~$k$ order
  isomorphic to the prefix (resp. suffix) of $\tau$ of length $k$. In
  case the prefix and suffix patterns of length $k$ of $\tau$ coincide, we say that $\tau$ has a \emph{bifix
  pattern of length} $k$.
\end{defn}

It is useful to note that in our poset $\mp$, each permutation can
cover at most two different permutations.  Namely, if $\sig$ is
covered by $\tau$ then $\sig$ clearly must occur as all but the first
or all but the last letter of $\tau$.  Thus, to obtain a permutation
covered by $\tau$ we can only remove the first or last letter of
$\tau$, thus yielding at most two different permutations.  Moreover,
the permutations obtained by removing the first and the last letter,
respectively, from $\tau$ are order isomorphic if and only if $\tau$
is \emph{monotone}, that is, if $\tau$ is either the increasing
permutation $123\dots n$ or the decreasing permutation
$n(n-1)\dots21$.

In the case where $\sig$ occurs precisely once in $\tau$, we show
that $\mst$ depends only on the lengths, $a$ and $b$, of the two
tails of $\tau$.  More precisely, $\mst$ is 1 if $a=b\le1$, it is
$-1$ if $a=0$ and $b=1$ or vice versa, and 0 otherwise (in which
case $\tau$ has a tail of length at least 2).

Our main result, Theorem~\ref{thm-x}, deals with intervals
$[\sig,\tau]$ where $\sig$ occurs at least twice in $\tau$.  This
result implies that, as in the case of one occurrence, if $\tau$
has a tail of length at least 2, then $\mst=0$.  In the remaining
cases, where the tails of $\tau$ have length at most 1, the main
result gives a recursive algorithm for computing $\mst$, by
producing, if possible, an element $\x$ in $\ist$, where
$|\x|<|\tau|-2$, such that $\mst=\mu(\sig,\x)$.  This element
$\x$, if it exists, must be a bifix pattern of $\tau$, and it must
lie below the two elements covered by $\tau$, but not below the
element obtained by deleting one letter from each end of $\tau$.
If no such element $\x$ exists (which is most often the case), we
have $\mst=0$.

\section{The case of one occurrence}\label{sec-one-occ}

First we consider the case when $\sigma$ occurs precisely once in
$\tau$.  In this case the interval $[\sigma,\tau]$ can be described
very simply.  Let $w_1\overline{\sigma}w_2$ be the factorization of
$\tau$, where the entries of $\overline{\sigma}$ constitute the only
occurrence of $\sig$ in $\tau$ and $|w_1|=a$, $|w_2|=b$, so $\tau$ has
left and right tails of lengths $a$ and $b$, respectively.  We obtain
the permutations covered by $\tau$ by deleting the first or the last
entry of $\tau$ and renaming the remaining entries.  Starting from
$\tau$ and deleting elements from both tails in such a way, we are
eventually left with $\sigma$.  It is easy to see that the interval
$(\sig,\tau)$ is a \emph{grid}, that is, a direct product of two
chains of lengths $a$ and $b$, respectively.  See
Figure~\ref{one_occurrence} for an example.  It is also easy to see
that the M\"obius function of $(\sig,\tau)$ for such a grid is 0,
except when $a$ and $b$ are both at most $1$.  This is recorded in the
following theorem.

\setlength{\unitlength}{5mm}
\begin {figure}[htbps]
\begin{picture}(20,13)(-2,-1)

\allinethickness{.2mm}

\def\xo{0}
\def\xa{3}
\def\xb{6}
\def\xc{9}
\def\xd{12}
\def\xe{15}

\def\yo{0}
\def\ya{2}
\def\yb{4}
\def\yc{6}
\def\yd{8}
\def\ye{10}

\def\cbox#1{\makebox(0,0){#1}}

\def\whitecbox#1#2#3#4#5{
\put(#1,#2){\cbox{\White{\rule{#3}{#4}}}}
\put(#1,#2){\cbox{#5}}
}

\def\wcbox#1#2#3{\whitecbox{#1}{#2}{20mm}{6mm}{#3}}

\path(\xo,\yb)(\xc,\ye)(\xe,\yc)(\xb,\yo)(\xo,\yb)
\path(\xa,\ya)(\xd,\yd)
\path(\xa,\yc)(\xc,\ya)
\path(\xb,\yd)(\xd,\yb)


\wcbox{\xc}{\ye}{68513427}

\wcbox{\xb}{\yd}{7513426}
\wcbox{\xd}{\yd}{6751342}

\wcbox{\xa}{\yc}{513426}

\wcbox{\xc}{\yc}{651342 }

\wcbox{\xe}{\yc}{564123}

\wcbox{\xo}{\yb}{13425}
\wcbox{\xb}{\yb}{51342}
\wcbox{\xd}{\yb}{54123}

\wcbox{\xa}{\ya}{1342}
\wcbox{\xc}{\ya}{4123}

\wcbox{\xb}{\yo}{123}

\end{picture}
\caption {The interval
  $[123, 68513427]$}\label{one_occurrence}
\end {figure}
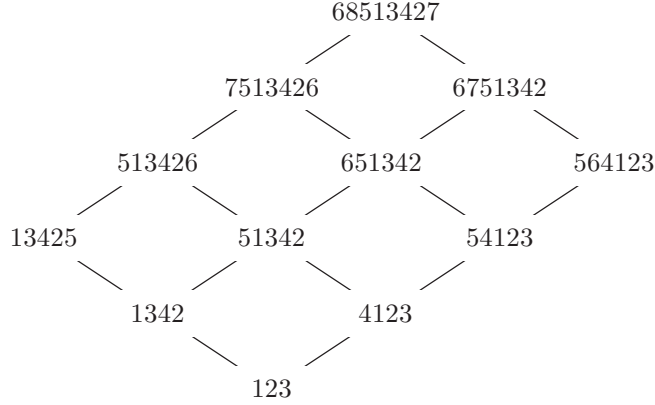

\begin{teor}\label{thm-one-occ}
Suppose $\sig$ occurs precisely once in $\tau$, and that $\tau$ has
tails of lengths $a$ and $b$, respectively.  Then $\mu(\sig,\tau)=1$
if $a=b=0$ or $a=b=1$, and $\mu(\sig,\tau)=-1$ if $a=0,b=1$, or
$a=1,b=0$.  Otherwise, $\mu(\sig,\tau)$ is 0.
\end{teor}

\section{More than one occurrence}

\begin{lemma}\label{lemma-short}
Suppose $\sig$ occurs in $\tau$ and that $r=|\tau|-|\sig|\le2$.  Then,
if $r=0$, $\mst=1$.  If $r=1$, we have $\mst=-1$.  If $r=2$, then
$\mst=0$ if $\tau$ (and thus also \sig) is a monotone permutation,
otherwise $\mst=1$.
\end{lemma}

\begin{proof}
If $r<2$, the claim is obvious, since $\ist$ is either a singleton or
a chain of two elements.  If $r=2$, and $\tau$ is monotone, removing a
letter from either end of $\tau$ yields the same (monotone)
permutation, so $\ist$ is a chain of three elements, and $\mst=0$.  If
$r=2$ and $\tau$ is not monotone, removing one letter from one end of
$\tau$ yields a permutation different from the one obtained by
removing a letter from the other end, for otherwise $\tau$ would be
monotone.  In that case, $\ist$ is a Boolean algebra with four
elements, whose M\"obius function is 1.
\end{proof}

Because of Lemma \ref{lemma-short}, we will from now on only consider
intervals $\ist$ of rank at least three, that is, where
$|\tau|-|\sig|\ge3$.  We also only need to consider pairs
$(\sig,\tau)$ such that $\sig$ occurs at least twice in $\tau$, since
the single-occurrence case is taken care of in Section
\ref{sec-one-occ}.

We will frequently refer to certain special elements of an interval
$\ist$ described in the following Definition.

\begin{defn}
  Given a permutation $\tau$, we let $\ptau$ be the standard form of
  $\tau$ after having removed its first letter, $\taup$ be the
  standard form of $\tau$ after having removed its last letter, and
  $\ptp$ be the standard form of $\tau$ after having removed both its
  first and last letter.  We refer to $\ptp$ as the \emph{interior of
    $\tau$}.
\end{defn}

Thus, for instance, if $\tau =68513427$, then $\ptau =7513426$,
$\taup =6751342$ and $\ptp =651342$.

\begin{lemma}\label{lemma-one-end}
Suppose $\sig$ occurs at least twice in $\tau$ and that
$|\tau|-|\sig|\ge3$.  If {} $\ptp$ does not lie in $\ist$, we have
$\mst=1$.
\end{lemma}

\begin{proof}
Observe first that $\tau$ cannot be monotone, since then it would be
possible to remove one letter from each end and get a permutation in
\ist, given that $|\tau|-|\sig|\ge3$.

Together with the hypothesis this implies that $\tau$ must have
precisely two occurrences of $\sigma$, which necessarily appear at its
two ends. Therefore $\ist$ consists of two chains (of equal lengths)
having in common only the minimum ($\sigma$) and the maximum
($\tau$). It is easy to see that the \mf\ of such an interval is 1.
\end{proof}

In the lemmas above, and in Section \ref{sec-one-occ}, we have
taken care of all intervals except those where $\tau$ has at least
two occurrences of $\sig$, $|\tau|-|\sig|\ge3$ and the interior of
$\tau$ lies in \ist.  We now deal with these remaining cases.

\begin{lemma}\label{lemma-C} Given $\sig$ and $\tau$ in $\mp$, let
$$ \pC=\{\rho\in[\sigma,\tau]\;|\;\rho<\ptau ,\rho\not\le\ptp\}$$ and
  let
$$ \Cp=\{\rho\in[\sigma,\tau]\;|\; \rho<\taup,\rho\not\le\ptp\}.$$ The
  sets $\pC$ and $\Cp$ are chains, and $\pC\cap\Cp$ has at most one
  element.  Moreover, if $z\in\Cp\setminus\pC$ or
  $z\in\pC\setminus\Cp$ then $[z,\tau]$ is a chain.
\end{lemma}

\begin{proof}
A permutation $\rho$ in $\pC$ cannot occur in the interior of
$\tau$ (since $\rho \not\le \ptp$), and thus it has to be a suffix
pattern of $\tau$ (since $\rho <\ptau$).  This implies that $\pC$
is a chain, since two suffix patterns of $\tau$ must be comparable
in $\mp$. The argument for $\Cp$ is analogous, so all permutations
in $\Cp$ occur as prefix patterns of $\tau$.

Suppose that $x,y\in \pC \cap \Cp$ and $x\neq y$.  We can assume,
without loss of generality, that $y<x$, since $x$ and $y$ are
elements of the chain $\Cp$.  Then $y$ must be a proper prefix
pattern of $x$, since both belong to $\Cp$ and thus are prefix
patterns of~$\tau$. Likewise, since both belong to $\pC$, $y$ must
be a proper suffix pattern of $x$.  But, an element that occurs as
a strict suffix of a prefix of $\tau$ and as a strict prefix of a
suffix of $\tau$ must occur in the interior of $\tau$, and thus we
must have $y\le\ptp$, which contradicts the assumption that
$y\in\Cp$.  So,  $\pC\cap\Cp$ can contain at most one element.

Assume now that $z\in\pC\setminus\Cp$. Let $y\in [z,\tau ]$ and
suppose that $y$ is neither $\ptau$ nor $\tau$ (in which cases
$[y,\tau]$ would be a chain). In this case $y$ cannot lie below
either $\taup$ or $\ptp$, since then we would have $z<\ptp$
(because $\pC$ is a chain), contradicting the hypothesis $z\in
\pC$. Hence, for all $y\in [z,\tau]$, we have
$y\in\pC\cup\{\ptau,\tau\}$, which is a chain. Then $[z,\tau]
\subseteq\pC\cup\{\ptau,\tau\}$, whence $[z,\tau]$ is a chain. An
analogous argument shows that $[z,\tau]$ is a chain if
$z\in\Cp\setminus\pC$.
\end{proof}

We have thus shown that $\pC \cap \Cp$ has at most one element.  In
case it exists, we give it a special name, and record its properties
in the following lemma, which is an immediate consequence of the proof
of Lemma~\ref{lemma-C}.

\begin{lemma}[and Definition]\label{lemma-x}
Given $\sig$ and $\tau$ in $\mp$, let $\pC$ and $\Cp$ be as defined in
Lemma~\ref{lemma-C}.  If {} $\pC\cap\Cp$ is nonempty, let $\x$ be its
(necessarily unique) element.  Then $\x$ has the following properties:
\begin{enumerate}
\item\label{item-1} $\x<\ptau$ and $\x<\taup$,

\item\label{item-2} $\x\not\le\ptp$.
\end{enumerate}

Conversely, if $\x$ satisfies the above two conditions, then $\x\in
\pC\cap\Cp$. In what follows, we refer to $\x$ as the \emph{carrier
  element of $\ist$}, or simply as $\x$, if it is clear from the
  context what $\sig$ and $\tau$ are.
\end{lemma}

Observe that $[\sigma ,\tau ]$ can be expressed as the disjoint
union of the principal ideal generated by $\ptp$ and the filter
$\pC \cup \Cp\cup\{\ptau,\taup,\tau\}$.  In case $\x$ exists, such
a filter is the principal filter generated by $\x$. This fact
holds in general, that is, even when the interval $\ist$ has rank
$r\le2$. More precisely, if $r=1$, then the ideal $I$ is empty and
the filter is $F=\{\sigma,\tau\}$, whereas if $r=2$, then
$I=\{\ptp\}=\{\sigma\}$ and $F=\{\ptau,\taup,\tau\}$.

\begin{teor}\label{thm-x}
  Suppose $\tau$ has at least two occurrences of $\sig$ and that
  $|\tau|-|\sig|\ge3$.  Assume that
  $\ptp$ 
  lies in \ist.  Then, if $\ist$ has no carrier element, we have
  $\mu(\sigma,\tau)=0$.  Otherwise, $\mu(\sigma,\tau)=\mu(\sigma,\x)$.
\end{teor}

\begin{proof}
  In our argument we will compute the M\"obius function of $\ist$ from
  top to bottom, using the formula for the M\"obius function in
  Equation~(\ref{up-down-mob}).  We will thus compute $\mst$ by
  recursively computing the value of $\m z,\tau,$ for $z$ of
  decreasing ranks, starting with $z=\tau$, as is done on the right
  hand side in Figure~\ref{fig-updown}.

  Note that, since $|\tau|-|\sig|\ge3$, we have that $\sig$ lies
  strictly below all of $\tau,\taup,\ptau,\ptp$, and $\sigma$ is in
  particular not equal to any one of them.

  If $z\in\Cp\setminus\pC$ then, by Lemma~\ref{lemma-C}, the interval
  $[z,\tau]$ is a chain.  Since $z<\taup<\tau$, this chain has rank at
  least 2.  Thus, $\mu (z,\tau )=0$.  An analogous argument shows that
  $\mu (z,\tau )=0$ when $z\in\pC\setminus\Cp$.  Thus, every element
  $t$ in $\Cp\cup\pC$, except $\x$ (if it exists), has
  $\mu(t,\tau)=0$.

  If $\x$ does not exist, then we claim that $\m y,\tau,=0$ whenever
  $y<\ptp$, that is whenever $y$ is different from $\ptp$ but belongs
  to the principal ideal generated by~$\ptp$.  We claim that a maximal
  such element $y$ (which must be covered by $\ptp$) lies below
  precisely four elements $t$ with $\mu(t,\tau)\neq0$, namely
  $\tau,\taup,\ptau,\ptp$, and hence has $\m y,\tau,=0$. This is
  because any other element $z$ that $y$ could lie below must belong
  to $\pC \cup \Cp$, since $z<\ptp$ is impossible by virtue of $y$
  being maximal.  Thus, $\mu(z,\tau)=0$, according to the preceding
  paragraph. By induction, this now also applies to all other elements
  lying strictly below $\ptp$, since each of them lies below precisely
  four elements $t$ with $\mu(t,\tau)\ne0$, namely, $\ptp,\ptau,\taup$
  and $\tau$. In particular, this shows that $\mst=0$ if $\ist$ has no
  carrier element, and, in fact, that then $\mu(z,\tau)=0$ for every
  $z\in\ist$ except for $\tau,\taup,\ptau$ and $\ptp$.

Finally, assume that $\x$ exists.  We claim that $\m \sigma, \tau,
=\m \sigma, \x,$.  Indeed, we have
\begin{displaymath}
\m \sigma, \tau, =-\sum_{\sigma <z\leq \tau}\m z, \tau, .
\end{displaymath}

If $z\not\le \x$, we have seen above that $\m z, \tau, =0$, whence
\begin{equation}\label{induction}
\m \sigma, \tau, =-\sum_{\sigma <z\leq \x}\m z, \tau, .
\end{equation}

Next we prove that $\m z,\tau,=\m z,\x,$ whenever $z\le \x$.  We
proceed by induction on the difference between the rank of $\x$ and
the rank of $z$.  If $z=\x$, then $\m z,\tau, =1$, since the only
elements between $\x$ and $\tau$ having nonzero \mf\ are $\tau ,\ptau$
and $\taup$.  This is the base case of the induction.  For any $z\leq
\x$, using the induction hypothesis on the elements between $z$ and
$\x$ and strictly above $z$, we have
\begin{displaymath}
\m z, \tau, =-\sum_{z<t\leq \tau}\m t, \tau, =-\sum_{z<t\leq \x}\m
t, \tau, =-\sum_{z<t\leq \x}\m t, \x, =\m z, \x, .
\end{displaymath}

Plugging this into formula (\ref{induction}) we get
\begin{displaymath}
\m \sigma, \tau, =-\sum_{\sigma <z\leq \x}\m z, \x, =\m \sigma,
\x, ,
\end{displaymath}
as desired.
\end{proof}

We thus get a recursion, where we find the carrier element of the
interval $[\sigma,\x]$, and iterate this until we get an interval that
does not have a carrier element.  That final interval either has
M\"obius function zero, or else its rank is at most 2, making it
trivial to compute the M\"obius function.

Concerning time complexity, a rough analysis shows that, in the
worst case, the recursive procedure described by the above theorem
is essentially bounded above by $|\tau |^3$. Indeed, supposing
that $|\sigma |=k$ and $|\tau |=n$, in order to compute the
possible carrier element of $\ist$, one first tries with the
prefix of $\tau$ of length $n-2$, and has to check both that it is
isomorphic to the suffix of $\tau$ of length $n-2$ and that it is
not isomorphic to $\ptp$. The time needed for this task is
proportional to $(n-2)^2$. In case no carrier element has been
found, one has to try with the prefix of length $n-3$; by an
analogous argument, the time needed is proportional to $(n-3)^2$.
In the worst case, (that is, if $\ist$ has no $\x$ and $\sigma$
occurs at both ends of $\tau$), this task has to be performed
until we get to $\sigma$ (i.e. the prefix of $\tau$ of length
$k$). So the total time needed in this worst case is proportional
to $\sum_{i=k}^{n-2}i^2$, which is proportional to~$n^3$.

Observe, however, that the worst case is not the most significant
one. Indeed, our results imply that it is difficult for the M\"obius
function of $\ist$ to be nonzero when $\sigma$ is very long. To be
more precise, when $\tau$ has a tail with respect to $\sigma$ of
length at least 2, then the M\"obius function of $\ist$ is equal to
0. Thus if $\sigma$ is not order isomorphic to four specific subwords
of $\tau$ (this is a quadratic test), then the M\"obius function is
0. Moreover, observe that the probability that $\sigma$ appears at the
beginning of $\tau$, for instance, is equal to $1/k!$ (since we have
to choose a $k$-subset of an $n$-set, then to arrange its elements in
the unique way which produces a word order-isomorphic to $\sigma$,
finally we can permute the remaining $n-k$ elements as we like). So,
when the length of the pattern $\sigma$ increases, the probability to
have $\sigma$ at the beginning of $\tau$ (as well as in any other
specific position) rapidly decreases.

In what follows we will frequently refer to the last $\x$ we find
in an interval $\ist$ when finding the carrier element $\x$ of
$\ist$, then the carrier element of $[\sigma,\x]$ and so on until
we have come to the last carrier element $\x'$ in this sequence
(which implies that $[\sigma,\x']$ does not have a carrier
element).  In this case, we will refer to $\x'$ as the
\emph{socle} of $\ist$.

\begin{cor}\label{long_tail} Suppose $\sig$ occurs in $\tau$ but that
the first two (or the last two) letters of $\tau$ are not involved in
any occurrence of \sig.  Then $\mst=0$.  \end{cor}

\begin{proof} If $\sig$ occurs only once in $\tau$ the claim has
  already been established, in Theorem \ref{thm-one-occ}.  Assume then
  that $\sig$ occurs at least twice in $\tau$.  Then
  $|\tau|-|\sig|\ge3$.  According to Theorem \ref{thm-x}, if $\ist$
  has no carrier element $\x$ satisfying the hypotheses of
  Theorem~\ref{thm-x}, then $\mst=0$.  Otherwise, if there is such a
  $\x$, then, as a consequence of its definition, $\x$ is order
  isomorphic to an initial segment of $\tau$, and to a final segment
  of $\tau$.  Thus, there is no occurrence of $\sig$ in $\x$ involving
  the first two (or the last two) letters of \x.  Iterate now the
  construction of $\x$ for $[\sig,\x]$ until we obtain the socle $\x'$
  of $\ist$. Since $\x'$ has a tail of size at least two with respect
  to $\sig$ we have $|\x'|-|\sig|\ge2$ and there are two
  possibilities: if $|\x'|-|\sig|=2$ it follows that $[\sig,\x']$ is a
  three element chain, so $\m\sig,\x',=0$. Otherwise,
  $|\x'|-|\sig|\ge3$, in which case Theorem \ref{thm-x} implies that
  $\m\sig,\x',=0$.
\end{proof}

Observe that if $\x$ exists then, since it lies in both $\pC$ and
$\Cp$, $\x$ must be a bifix pattern of $\tau$ of length at least
$a+|\sigma |+b$, where $a$ and $b$ are the lengths of the left and
right tails (respectively) of $\tau$ with respect to $\sigma$.
Also, if $\x$ exists, it does not appear anywhere else in $\tau$
(this just follows from the definition of $\tau$). Here are two
examples of the $\x$ associated to an interval: For $\sig=321$ and
$\tau=431825976$, we have $\x=321=\sig$, so $\mst=\m\sig,\sig,=1$,
whereas if $\sig=231$ and $\tau=2\;5\;7\;1\;4\;8\;9\;3\;6\;10$ we
have $\x=245136$ and $\mst=\msx=0$ since the interval
$[231,245136]$ has no carrier element. Observe that in the latter
case the initial and final segments of $\tau$ order isomorphic to
$\x$ overlap in the letters 48.

\bigskip

As a consequence of Corollary~\ref{long_tail}, nonzero values of
$\mst$ can occur only when the tails of $\tau$ have length at most
1. Below we will always assume that $\tau$ satisfies this hypothesis,
and we will give a series of partial conditions that facilitate the
computation of the M\"obius function in several cases.

Denote with $x$ the sum of the lengths of the tails of $\tau$ with
respect to $\sigma$. So $x=0$ means that $\tau$ has two occurrences of
$\sigma$, one at each end; $x=1$ means that $\tau$ has one occurrence
of $\sigma$ at one end and a tail of length 1 at the other end;
finally, $x=2$ means that $\tau$ has two tails of length 1 each.

For simplicity, in what follows we will always assume that, when
$x=1$, $\tau$ has an occurrence of $\sigma$ at its right end (and
thus a tail of length 1 at its left end). In case $\sigma$ appears
at the left end of $\tau$, we simply have to replace each
occurrence of the word ``suffix" with the word ``prefix" in all
the following propositions.

Our first result says that it is very difficult for the M\"obius
function to take the value $(-1)^{x+1}$.

\begin{prop}\label{first}
\begin{enumerate}
\item[1.] Suppose $x=0$. If $\tau$ does not have a monotone bifix of
  length $|\sigma |+1$, then $\mst \neq -1$.
\item[2.] Suppose $x=1$. If $\tau$ does not have a bifix of length
  $|\sigma |+2$ whose suffix of length $|\sigma |+1$ is monotone, then
  $\mst \neq 1$.
\item[3.] Suppose $x=2$. Then $\mst \neq -1$.
\end{enumerate}
\end{prop}

\begin{proof}
  Observe that, in general, any permutation of length $\ell$ having a
  bifix of length $\ell-1$ is necessarily monotone.  Suppose $x=0$. If
  $\mst =-1$, then (by Theorem \ref{thm-x} and Lemma
  \ref{lemma-short}) if $\ist$ has a carrier element, then the socle
  of $\ist$ must have length $|\sigma|+1$, so it is necessarily
  monotone (since $\sigma$ is a bifix of it). Suppose $x=1$. If $\mst
  =1$, then (by virtue of Theorem \ref{thm-x} as well) the socle must
  have length $|\sigma |+2$, and a direct inspection shows that its
  suffix of length $|\sigma|+1$ has to be monotone (since $\sigma$ is
  a bifix of it). Suppose $x=2$. Then necessarily the socle (if it
  exists) has length at least $|\sigma|+2$, and so $\mst \neq -1$.
\end{proof}

In this direction, a general result that includes all the cases of
(but is weaker than) the previous proposition is the following. The
proof is easy and is omitted.

\begin{prop} If $\tau$ has a non-monotone suffix of length $|\sigma
  |+x$, then $\mst \neq (-1)^{x+1}$.
\end{prop}

Next we give an easy necessary condition in order to have $\mst
=0$. For this, we first need a definition. A permutation is said
to be \emph{monotone (reverse) alternating} when it is (reverse)
alternating and the two permutations induced by its even-indexed
elements and odd-indexed elements are both monotone. For instance,
the permutation 342516 is monotone alternating. Alternating
permutations have been extensively studied also in recent years,
see for instance the survey \cite{S}.

\begin{prop}
\begin{enumerate}
\item[1.] Suppose $x=0$ and that $\sigma$ is not the socle of
  $\ist$. If~$\tau$ has neither a monotone bifix of length $|\sigma
  |+1$, nor a monotone (reverse) alternating bifix of length $|\sigma
  |+2$, then $\mst =0$.
\item[2.] Suppose $x=1$. If $\tau$ has neither a bifix of length
  $|\sigma |+1$ nor a bifix of length $|\sigma |+2$ whose suffix of
  length $|\sigma |+1$ is monotone, then $\mst =0$.
  \item[3.] Suppose $x=2$. If $\tau$ does not have a bifix of length
    $\sigma |+2$, then $\mst =0$.
\end{enumerate}
\end{prop}

\begin{proof}
  Suppose that $\mst \neq 0$. Then $\ist$ necessarily has a carrier
  element, and we let $\x'$ be the socle of $\ist$. It is clear that,
  in order to have $\mst \neq 0$, $\x'$ must be ``short'', and more
  precisely $|\x'|\leq |\sigma|+2$.

  If $x=0$, we have two distinct cases (since, by hypothesis,
  $|\x'|\neq |\sigma |$). 
  If $|\x'|=|\sigma |+1$, then $\sigma$ is a bifix of $\x'$, and so
  necessarily $\x'$ is a monotone bifix of $\tau$. Finally, assume
  that $|\x'|=|\sigma |+2$. Of course, also in this case $\sigma$ is a
  bifix of $\x'$. Suppose $\sigma =a_1 \cdots a_k$. Our hypothesis
  implies that, for all $i$, (i) $a_{2i}<a_{2i+1}$ if and only if
  $a_{2i+2}<a_{2i+3}$, and (ii) $a_{2i+1}<a_{2i+2}$ if and only if
  $a_{2i+3}<a_{2i+4}$. Therefore the descents of $\sigma$ are
  completely determined by the first three elements $a_1 ,a_2$ and
  $a_3$. Moreover, observe that, if we had either $a_1 <a_2 <a_3$ or
  $a_1 >a_2 >a_3$, then $\x'$ (and a fortiori $\sigma$) would be
  monotone, and this contradicts the fact that $[ \sigma ,\x']$ has no
  carrier element. Thus we must have either $a_1<a_2>a_3$ or
  $a_1>a_2<a_3$, which implies that $\x'$ is alternating. Finally,
  there cannot exist any $i$ such that $a_{2i}<a_{2i+2}$ and
  $a_{2i+2}>a_{2i+4}$, since the word $a_1 \cdots a_{k-2}$ is order
  isomorphic to the word $a_3 \cdots a_k$.  The remaining cases can be
  treated analogously.

If $x=1$, we have only two possibilities. If $|\x'|=|\sigma |+1$,
then $\x'$ is clearly a bifix of $\tau$ of length $|\sigma |+1$. If
$|\x'|=|\sigma |+2$, then Proposition \ref{first} tells us that
$\x'$ is a bifix of $\tau$ whose suffix of length $|\sigma |+1$ is
monotone.

If $x=2$, then necessarily $|\x'|=|\sigma |+2$, and $\x'$ is of
course a bifix of $\tau$.
\end{proof}

A general result in this direction is the following.

\begin{prop}
  Denote with $\omega$ the longest bifix of $\tau$ containing $\sigma$
  and having length $\leq |\sigma |+2$. Moreover, denote with $\alpha$
  (resp. $\beta$) the shortest prefix (resp. suffix) pattern of $\tau$
  containing the first (resp. last) two occurrences of $\omega$. If
  $\omega$ exists and $\alpha \neq \beta$, then $\mst =0$.
\end{prop}

\begin{proof}(Sketch) Using an argument similar to that of the
  previous proposition, suppose that $\mst \neq 0$, and let $\x''$ be
  the element of $\ist$ such that the socle of $\ist$ is the carrier
  element of the interval $[\sigma ,\x'']$. (Thus, $\x''$ is the
  penultimate carrier element found in the iteration culminating in
  the socle of $\ist$).  Clearly $\x''$ has to be a bifix of $\tau$,
  and it is easy to see that $\x''=\alpha =\beta$.
\end{proof}

\end{document}